\documentclass[11pt]{article} 

\usepackage{latexsym}
\usepackage{amsmath}
\usepackage{amssymb}
\usepackage{stmaryrd}

\newcommand{\vs}[1]{\vspace{#1mm}}
\newcommand{\hs}[1]{\hspace*{#1mm}}

\begin{document}

\begin{center}

{\bf \Large The diagonal lemma as the formalized Grelling 
paradox}\,\footnote{To appear in the series Collegium Logicum of the 
Kurt G\"odel Society in 2006.}

\vs{2}

{\bf Gy\"orgy Ser\'eny}

\end{center}

\vs{2}

\hs{3}G\"odel's diagonal lemma (which is often referred to 
as fix-point or self-referential lemma) summarizes very succinctly 
the ability of first-order arithmetic to `talk about 
itself', a crucial property of this system 
that plays a key role in the proof of 
G\"odel's incompleteness theorem and in those of the 
theorems of Tarski and Church on the undefinability of 
truth and undecidability of provability respectively. 
In fact, with the representability of provability at hand, 
these three main limitative theorems of logic 
can be considered to be simple applications of the lemma 
(cf.\,e.g.\,[2]\,pp.\,227--231). 
Due to this central role, the proof of the lemma could shed light 
on the very essence of these fundamental theorems. 
In spite of the fact that it is common knowledge that 
G\"odel's proof of the incompleteness theorem 
is closely related to the Liar paradox, 
the proof of the lemma as it is presented in textbooks
on logic is not self--evident to say the least. Indeed, 
in the {\it Handbook of Proof Theory}, the proof of 
the lemma is introduced by the following remark
(see [1], p.119): `This proof [is] quite simple but rather 
tricky and difficult to conceptualize.' Or to quote another opinion, 
`The brevity of the proof does not make for 
transparency; it has the aura of a magician's trick' 
(cf.\,[4],\,p.\,1). It seems, therefore, that the words of a respected 
logician reflect a widespread attitude to the proof of the lemma 
(see [5]): `[This] result is a cornerstone of modern logic. 
[...] You would hope that such a deep theorem would have an 
insightful proof. No such luck. [...] 
I don't know anyone who thinks he has a fully satisfying understanding 
of why the Self-referential Lemma works. It has a 
rabbit-out-of-a-hat quality for everyone.' 

\hs{3}In view of these remarks, we think that it is worth 
drawing attention to a possibility of making the proof of the
lemma completely transparent by showing that it is simply a 
straightforward translation of 
the Grelling paradox into first-order arithmetic.\footnote{The 
train of thought below is an application of the ideas given in [6] 
to first-order arithmetic.} 

\vs{2}

{\bf Notation}

Our formal language is that of first-order arithmetic. $Q$ 
stands for Robinson arithmetic while $\omega$ is the set 
of natural numbers. $g$ is any one of the standard 
G\"odel numberings and $Fm_n$ is the set of formulas with 
all free variables among the first $n$ ones.
For the sake of simplicity, we shall denote the closed 
terms corresponding to natural numbers
by the numbers themselves. Further,  $N$ denotes the set of 
G\"odel numbers of formulas in $Fm_1$. Finally, 
the result of substituting a term $t$ 
for the only free variable of a formula $\varphi\in Fm_1$ 
is denoted by $\varphi(t)$.

\vs{2}

{\bf Diagonal lemma}  

\vs{1}

For any formula $\varphi\in Fm_1$, there is a 
sentence $\lambda$ such that 

\vs{2}

\ \hfill $Q\vdash \lambda \longleftrightarrow
\varphi(g(\lambda))\,.$ \hfill \ 

{\bf Proof idea} 

\vs{1}

First we show how to construct, out of Grelling's paradox, 
an ordinary language sentence that, on the one hand, 
says of itself that it has a given property, on the other hand, 
consists of components with 
easily identifiable formal first-order counterparts.
The straightforward formalization of this ordinary language 
sentence leads to the desired formal sentence 
(as can be expected since the lemma is just about the existence of 
a first-order sentence that, informally speaking, says of itself that 
it has a given property).

\hs{3}As is well known, the Grelling paradox consists in  
the fact that the sentence 

\vs{1}

(1) \hfill {\small `heterological' is 
heterological}\footnote{An adjective is called {\it heterological}
if the property denoted by the adjective 
does not hold for the adjective itself; e.g. 
`long', `German', `monosyllabic' are 
heterological.} \hfill \ 

\vs{1}

shares with the Liar sentence the remarkable property that 
its truth implies its own falsity and {\it vice versa},
i.e., in effect, says of itself that it 
is false. What is truly important is that, contrary to the Liar,  
this paradoxical sentence achieves self-reference 
without using an indexical.\footnote{As to the notion of 
heterologicality itself (which actually involves self-reference), 
its slightly modified version can also be 
expressed without using an indexical: replace `it' by `$x$' 
in (3) below.} 

\hs{3}Since our aim is to construct a sentence that
(a) is not about an adjective but about a sentence and
(b) instead of asserting its own falsehood, 
says of itself that it has an arbitrary (but 
fixed) property, 
we have to slightly modify (1) accordingly. In order to satisfy 
the first requirement, in place of an adjective $A$, we consider 
the open sentence `$x$ is $A$'. Obviously, in this case,  
the transformation corresponding to the application 
of an adjective $A$ to a linguistic object $O$ will be the 
substitution of the name of 
$O$\,\footnote{Following the common practice, we define 
the name of a linguistic object to be the object 
itself between quotation marks.} 
for the variable $x$ in the open sentence 
corresponding to $A$. Consequently, the sentence 
associated with the self-application 
of any adjective $A$ in this way is 
``\,`$x$ is $A$' is $A$''. In particular, 
the counterpart of (1) is:

\vs{1}

(2) \hfill {\small `$x$ is heterological' is heterological.}\hfill \ 

\vs{1}

Note that the notion of heterologicality occurring here is
already a property of open sentences with single variables.
Since, on the one hand, to be heterological is to 
have the property that its application to itself 
yields a false sentence, on the other, as we noted above, 
in the case of sentences, `applied to itself'  means 
`its name is substituted for the variable in it', 
for any open sentence $x$ with a single variable, we have

\vs{1}

(3) \hfill \begin{minipage}[t]{13cm}
{\small $x$ is heterological just in case 
the sentence obtained by 
substituting the name of 
$x$  for the variable in it is false.}
\end{minipage} \hfill \ 

\vs{2}

Finally, if we replace `being false' by `having property $p$', 
(2) and (3) together  yield:

\vs{0.5}

(4)\hfill 
\begin{minipage}[t]{13cm}  
{\small the sentence obtained by substituting the name of 
\,\,`the sentence obtained by 
substituting the name of 
$x$  for the variable in it has property $p$'\,\, 
for the variable in it has property\,\,$p$.}
\end{minipage} \hfill \

\vs{1.5}

It can directly be checked that this sentence indeed says of 
itself that it has property $p$ 
(and says nothing else) since it is built up in such a way that if we 
perform the substitution described in it, then we obtain the sentence 
itself, which is stated to have  
property\,\,$p$.\footnote{J.N. Findlay used sentences of 
the same structure to examine 
  informally the incompleteness 
theorem (cf.\,[3]).}  
Now, let $s$ denote the open sentence between the quotation 
marks in (4):

\vs{0.5}

\ \hfill {\small the sentence obtained by substituting the name of 
$x$  for the variable in it has property $p$.} \hfill \ 

\vs{0.5}

Then, clearly, the whole sentence (4) is  
$s(`s\mbox{'})$.\footnote{Mimicking the 
formal notation, in the case of any common language 
open sentence $o$ having a single 
variable, we abbreviate the result 
of substituting a linguistic phrase  $q$  for the variable 
in $o$ by $o(q)$.} 
That is, the formalization process should consist of two steps. 
In the first step we have to find the formal version $\eta\in Fm_1$ 
of $s$, and then the second step is obvious: the desired sentence 
$\lambda$ will simply be $\eta(g(\eta))$.\footnote{G\"odel
numbering is, of course, the formal counterpart of naming.}

\vs{2}

{\bf Proof}

\vs{1}

Let $\varphi\in Fm_1$ be arbitrary and let its informal 
counterpart be the open sentence `$x$ has 
property $p$'.\footnote{It is obvious that the 
formulas in $Fm_1$ are formal versions of open sentences 
with single variables asserting the possession of a property, 
and, taking into consideration only those 
informal concepts that have formal counterparts, 
the formalization of attributing a property to an object is 
the substitution of the formal name (i.e. the G\"odel number) 
of the corresponding formal object for the only 
free variable of the formula that formalizes the open sentence 
asserting the possession of the property concerned.}
Certainly, $x(g\mbox{{\small (}}x\mbox{{\small )}})$ 
is the formal version of the phrase 

\vs{0.3}

\ \hfill {\small the sentence obtained by 
substituting the name of $x$  for the variable in it,} \hfill \ 

\vs{0.5}

and hence $\varphi\big(g[x(g\mbox{\small (}x\mbox{\small )})]\big)$ 
is the formal version of $s$. 
Clearly, $\varphi\big(g[x(g\mbox{\small (}x\mbox{\small )})]\big)$ 
with a variable $x$ running over formulas in 
$Fm_1$\footnote{Recall that, in $s$, 
the variable $x$ runs over the set of open sentences 
with single \vs{-2}variables.}
is not a formula itself, it becomes a formula only if we
replace the variable $x$ by a formula. Therefore, 
we cannot continue the formalization process unless
we find a formula that can play the role of 
$\varphi\big(g[x(g\mbox{\small (}x\mbox{\small )})]\big)$, 
that is, a formula $\eta\in Fm_1$ such that 
$\eta(g(\psi))$ is provably equivalent in $Q$ to 
$\varphi\big(g[\psi(g\mbox{\small (}\psi\mbox{\small )})]\big)$ 
{\it for every $\psi \in Fm_1$},  
or equivalently (denoting the inverse of $g$ by $g^{-1}$), 
{\it for any $n \in N$},

\vs{1}

\ \hfill $Q\vdash  
\eta(n)\longleftrightarrow\varphi\big(g[g^{-1}(n)(n)]\big)\,. $
\hfill \ 

\vs{1}

In order to find the appropriate formula $\eta$, 
let us consider the expression
substituted into the formula $\varphi$, 
and define the function $f:\omega\longrightarrow \omega$ accordingly:

\vs{1}

\ \hfill $f(n)=g[g^{-1}(n)(n)]$ \ if \ $n\in N$ \ \ and \ \ $f(n)=0$ \ 
otherwise. \hfill \ 

\vs{1}

Since this function is obviously recursive and hence representable 
in $Q$, and, up to provable equivalence in $Q$, the result of 
substituting a representable function into a formula can also 
be expressed by a formula,\footnote{By elementary first-order logic, 
it follows from the definition of representability that if a function 
$f:\omega\longrightarrow \omega$ is 
represented  in $Q$ by a formula $\mu\in Fm_2$, then, 
for any $\varphi\in Fm_1$ and $n\in \omega$, 

\ \hfill $Q\vdash  
(\exists y)(\mu(n,y)\land \varphi(y))
\longleftrightarrow\varphi(f(n))$.\hfill \ }
 there is a formula $\eta\in Fm_1$
such that, for any $n\in N$, 

\vs{1}

(5) \hfill $Q\vdash  \eta(n)\longleftrightarrow\varphi(f(n))$\,. 
\hfill \ 

\vs{1}

Thus we have obtained what we need, we have shown that there exists an 
$\eta\in Fm_1$ that can be considered to be the formal version 
of $s$. Now, all that remains to do is straightforward: 
it follows from (5) that, for every $\psi \in Fm_1$,

\vs{1}

\ \hfill $Q\vdash \eta(g(\psi))\longleftrightarrow
\varphi\big(g[\psi(g\mbox{\small (}\psi\mbox{\small )})]\big)$,  
\hfill \ 

\vs{2}

which, in turn, choosing $\psi$ to be $\eta$, yields

\vs{2}

\ \hfill $Q\vdash \eta(g(\eta))\longleftrightarrow
\varphi\big(g[\eta(g\mbox{\small (}\eta\mbox{\small )})]\big)$,  
 \hfill \ 

\vs{2}

showing that the sentence 
$\lambda =\eta(g(\eta))$ indeed  has the desired 
property.\footnote{Perhaps it is worth noting that
the informal version of the last step in the formal proof
explains the reasons why $s$ (the informal counterpart of $\eta$) 
is suitable for constructing 
the appropriate self-referring sentence (4). Actually, by definition, 
for {\it any} open sentence $o$ with a single 
variable, $s(`o\mbox{'})$ says that $o(`o\mbox{'})$ 
has property $p$. In the particular case when $o$ is just $s$,
we obtain: \,$s(`s\mbox{'})$ says 
that $s(`s\mbox{'})$ has property $p$.}

\vs{6}

{\bf Acknowledgments}

\vs{1}

This work was supported by Hungarian NSF grant No. T43242.

\vs{4}

{\bf References}

{\small 

\vs{2}

[1]\hs{2}S.\,R.\,Buss (ed.), {\it Handbook of Proof theory}, 
Elsevier, Amsterdam, 1998.

[2]\hs{2}H.\,B.\,Enderton, {\it A Mathematical Introduction to Logic}, 
Academic Press, New York, 1972.

[3]\hs{2}J.\,N.\,Findlay, Goedelian Sentences: A Non-numerical 
Approach, {\it Mind}, Vol.\,51, 1942, pp.\,259-65.

[4]\hs{2}H.\,Gaifman, Naming and Diagonalization, 
To appear in the {\it Logic Journal of the IGPL}, 22 pp.,
{\tt http://www.columbia.edu/{\small $\sim$}hg17/naming-diag.pdf}

[5]\hs{2}V.\,McGee, G\"odel's First Incompleteness Theorem,  
Handout for the course 24.242 Logic II, Spring 2002, 
{\tt http://web.mit.edu/24.242/www/1stincompleteness.pdf}

[6]\hs{2}G.\,Ser\'eny, \,\,G\"odel, Tarski, Church, and the Liar, 
{\it The Bulletin of Symbolic Logic}, Vol\,\,9, 2003, pp.\,3-25.}

\vs{6}

{\small\it Department of Algebra

Budapest University of Technology and Economics

1111 Stoczek u. 2. H \'ep. 5. em., Budapest, Hungary

\vs{-0.5}

e-mail: {\sf sereny@math.bme.hu}}

\end{document}